\documentstyle[12pt,amssymb,draft]{amsart}
\theoremstyle{plain}
\newtheorem{thm}{Theorem}[section]
\newtheorem{theorem}[thm]{Theorem}
\newtheorem{corollary}[thm]{Corollary}
\newtheorem{lemma}[thm]{Lemma}
\newtheorem{remark}[thm]{Remark}

\newenvironment{theorem*}[1]{\smallskip\noindent{\bf #1.}\it}{\medskip}
\newenvironment{proofof}[1]{\smallskip\noindent{\it #1}\rm}
                {\hspace*{\fill} $\Box$\medskip}
\newenvironment{proof}{\smallskip\noindent{\it Proof.}\rm}
                        {\hspace*{\fill} $\Box$\medskip}

\numberwithin{equation}{section}
\newcommand\dd{\frac{d^2}{dt^2}}
\newcommand\ov{\overline}
\newcommand\ti{\widetilde}
\newcommand\lan{\langle}
\newcommand\ran{\rangle}
\newcommand\bn{\binom}
\newcommand\tr{\operatorname{tr}}

\renewcommand\Re{\operatorname{Re}}

\newcommand\al{\alpha}
\newcommand\be{\beta}
\newcommand\ga{\gamma}
\newcommand\de{\delta}

\newcommand\eps{\varepsilon}
\newcommand\la{\lambda}
\newcommand\si{\sigma}
\newcommand\th{\theta}

\newcommand\bC{{\mathbb C}}
\newcommand\bR{{\mathbb R}}
\newcommand\bZ{{\mathbb Z}}

\newcommand\fD{{\frak D}}

\newcommand\op{operator}
\newcommand\sa{selfadjoint}
\newcommand\ev{eigenvalue}

\title{Schr\"odinger operators with singular Gordon potentials${}^{\dag}$}
\thanks{${}^{\dag}$The work was partially supported by Ukrainian State Foundation
for Basic Research}

\author{R.~O.~Hryniv \and Ya.~V.~Mykytyuk}

\address{Institute for Applied Problems of Mechanics and Mathematics,
	3b~Naukova str., 79601 Lviv, Ukraine}
\email{hryniv@@mebm.lviv.ua}

\address{Department of Mechanics and Mathematics, Lviv National University,
	1 Universytetska str., 79602 Lviv, Ukraine}
\email{mykytyuk@@email.lviv.ua}

\date{\today}

\begin{document}

\begin{abstract}
Singular Gordon potentials are defined to be distributions from the space
$W^{-1}_{2,unif}(\bR)$ that are sufficiently fast approximated by periodic
ones. We prove that Schr\"odinger operators with singular Gordon potentials
have no point spectrum and show that a rich class of quasiperiodic
distributions consists of singular Gordon potentials.
\end{abstract}

\maketitle

\section{Introduction}\label{sec:intr}

In the Hilbert space $L_2(\bR)$, we consider a Schr\"odinger operator
\begin{equation}\label{eq:S}
	S = - \dd + q
\end{equation}
with potential $q$ that is a real-valued distribution from the space
$W^{-1}_{2,unif}(\bR)$. The exact definition of this space is given in
Section~\ref{sec:pre} below, and at this point we remark that any
$q \in W^{-1}_{2,unif}(\bR)$ can be represented (not uniquely) in the form
$q = \si' + \tau$, where $\si$ and $\tau$ are real-valued functions from
$L_{2,unif}(\bR)$ and  $L_{1,unif}(\bR)$, respectively, i.e.,
\begin{align*}
  \|\si\|^2_{2,unif} &:= \sup_{t\in\bR} \int_t^{t+1} |\si(s)|^2 ds < \infty,\\
   \|\tau\|_{1,unif} &:= \sup_{t\in\bR} \int_t^{t+1} |\tau(s)|\,ds < \infty,
\end{align*}
and the derivative is understood in the sense of distributions.
For singular potentials of this type the domain and the action of the \op~$S$
are given by
\[
	\fD(S) = \{f \in W_2^1(\bR) \mid
		f^{[1]}:= f'-\si f \in W^1_{1,loc}(\bR),\
			l(f) \in L_2(\bR)\}
\]
and
\[
	S f = l (f),
\]
where
\[
	l (f) := - (f' -\si f)' - \si f' + \tau f.
\]
It is easily seen that $l(f) = -f'' + qf$ in the sense of
distributions, which implies, firstly, that the operator~$S$ so defined
does not depend on the particular choice of $\si \in L_{2,unif}(\bR)$ and
$\tau\in L_{1,unif}(\bR)$ in the decomposition $q=\si'+\tau$ and,
secondly, that for regular potentials $q \in L_{1,loc}(\bR)$ the above
definition coincides with the classical one.

The above regularisation by quasi-derivative procedure was first suggested
in~\cite{AEZ} for the potential $1/x$ on $[-1,1]$ and then systematically
developed in~\cite{SS} for a general class of singular (i.e., not necessarily
locally integrable) potentials from $W^{-1}_{2,loc}(\bR)$. We observe that
this class includes, among others, the classical examples of Dirac and Coulomb
potentials that are used to model zero- and short-range interactions in
quantum mechanics. On the other hand, this approach is not applicable to
locally more singular interactions like $\de'$ or $1/|x|^{\ga}$ with
$\ga\ge3/2$; see, e.g., \cite{AGHH} and the references therein and in
\cite{SS,HM} for some of the related literature on Schr\"odinger operators
with singular potentials.

Potentials $q$ from the space $W^{-1}_{2,unif}(\bR)$ were considered in detail
in~\cite{HM}. In particular, it was shown there that the corresponding
\op~$S$ is \sa, bounded below, and coincides with the form-sum of
the operator $-\dd$ and the multiplication operator by~$q$.
In the present paper we shall be interested in spectral properties of
the \op~$S$; more precisely, we would like to find conditions on the potential
$q \in W^{-1}_{2,unif}(\bR)$ under which the \op~$S$ does not have point
spectrum.

It is a classical result, for instance, that the spectrum of $S$
is purely absolutely continuous if the potential $q$ belongs to
$L_{1,loc}(\bR)$ and is periodic. Gordon~\cite{G} introduced a class of
bounded below $L_{2,loc}(\bR)$-potentials that are well enough approximated
by periodic ones and for which the corresponding Schr\"odinger \op s do not
possess eigenvalues. Damanik and Stolz~\cite{DS} enlarged recently this
class to $L_{1,loc}(\bR)$-valued functions. More exactly, they defined $q$ to
be a {\em generalized Gordon potential} if $q \in L_{1,unif}(\bR)$ and there
exist $T_m$-periodic functions $q_m\in L_{1,loc}(\bR)$ such that
$T_m \to \infty$ and for an arbitrary $C>0$ it holds
\[
	\lim_{m\to\infty} \exp(CT_m) \int_{-T_m}^{2T_m} |q(t) - q_m(t)|\,dt=0.
\]
It is proved in~\cite{DS} that the Schr\"odinger \op s with generalized
Gordon potentials have empty point spectrum. In particular, the authors
studied potentials of one-dimensional quasicrystal theory~\cite{AGHH} of the
form
\begin{equation}\label{eq:q-per}
	q(t) = q_1(t) + q_2(\al t + \th),
\end{equation}
where the functions $q_1,q_2\in L_{1,loc}(\bR)$ are $1$-periodic and
$\al,\th \in [0,1)$, and showed that~\eqref{eq:q-per} is a generalized Gordon
potential if $\al$ is a Liouville number~\cite{Lio} and $q_2$ is either
(i) H\"older continuous or (ii) a step function or (iii) $|x|^{-\ga}$ on
$[-1/2,1/2]$ with $\ga\in(0,1)$ or (iv) a linear combination of such
functions.

The class of generalized Gordon potentials, even though being large enough,
does not contain, for example, potentials modelling $\de$-interactions
that arise in one-dimensional crystal and quasicrystal theory
(cf. Kronig-Penney models and their various generalizations
in~\cite[Ch.III.2]{AGHH}, \cite{GK}, \cite{BSW}). The main aim of the present
article is to describe singular potentials in~$W^{-1}_{2,unif}(\bR)$
(with say $\de$-like singularities)
for which the above-defined operator~$S$ has empty point spectrum.
As it is shown in~\cite{HM}, so are all periodic potentials
in~$W^{-1}_{2,unif}(\bR)$; here we introduce singular analogues of Gordon
potentials that possess the property mentioned.

Namely, we say that $q = \si' + \tau \in W^{-1}_{2,unif}(\bR)$
with $\si\in L_{2,unif}(\bR)$ and $\tau \in L_{1,unif}(\bR)$
is a {\em singular Gordon potential} if there exist $T_m$-periodic
functions $\si_m \in L_{2,unif}(\bR)$ and $\tau_m \in L_{1,unif}(\bR)$
such that $T_m \to \infty$ and
\[
   \lim_{m\to\infty}\ \exp(CT_m) \Bigl\{ \|\si - \si_m\|_{L_2[-T_m, 2T_m]}
	+\|\tau - \tau_m\|_{L_1[-T_m, 2T_m]} \Bigr\} = 0
\]
for every $C < \infty$. (In fact, Theorem~\ref{thm:main} below remains true
if we only require this relation to hold for all $C < C_q$, where $C_q$ is
some constant dependent on $q$; see Remark~\ref{rem:q}.)
Observe that every generalized Gordon potential is evidently
a singular Gordon one and that $q = \si'+ \tau$ is a singular Gordon
potential if and only if
\[
   \liminf_{T\to\infty}\frac1T\log \Bigl\{ \|\si - \si_{(T)}\|_{L_2[-T, T]}
	+\|\tau - \tau_{(T)}\|_{L_1[-T, T]} \Bigr\} = - \infty;
\]
here $\si_{(T)}(t) := \si (t+T)$ and $\tau_{(T)}(t):=\tau(t+T)$.

The main results of the present paper are as follows.

\begin{theorem}\label{thm:main}
Suppose that $q$ is a singular Gordon potential. Then the \op~$S$
in~\eqref{eq:S} does not have any \ev s.
\end{theorem}

A rich class of singular Gordon potentials (in addition to those of the
form~\eqref{eq:q-per} pointed out in the work~\cite{DS}) is formed by
functions
\begin{equation}\label{eq:qper}
    q(t) = \si'_1(t) + \si'_2(\al t + \th) + \tau_1(t) + \tau_2(\al t + \th),
\end{equation}
where $\si_1,\si_2 \in L_{2,loc}(\bR)$ and $\tau_1, \tau_2\in L_{1,loc}(\bR)$
are $1$-periodic and $\al,\th \in [0,1)$.
Namely, we prove the following statement.

\begin{theorem}\label{thm:sgp}
Under the above assumptions the function $q$ of~\eqref{eq:qper}
is a singular Gordon potential if $\al$ is a Liouville number and $\si_2$
belongs to the space $W^s_{2,loc}(\bR)$ with some $s>0$.
\end{theorem}

In particular, Theorem~\ref{thm:sgp} shows that the requirement of~\cite{DS} that
the function~$q_2$ of~\eqref{eq:q-per} be of the above form (i)--(iv)
for $S$ to have empty point spectrum is superfluous; in addition, it
may have power-like singularities $|x|^{-\ga}$ with $\ga\in(0,3/2)$.

Absence of the point spectrum for the \op~$S$ is proved by the scheme analogous
to the one used in~\cite{G,DS}. Since $q=\si'+\tau$ is well approximated by
periodic functions $q_m = \si_m' + \tau_m$, a Gronwall-type inequality in
conjunction with some growth estimate gives proximity of solutions of the
equations
\begin{align}\label{eq:nprd}
	- (u' - \si u)' - \si u' + \tau u &= \la u,\\
	- (u_m' - \si_m u_m)' - \si_m u_m + \tau_m u_m &= \la u_m,\label{eq:prd}
\end{align}
satisfying the same initial data. If $u$ is an eigenfunction of the \op~$S$,
then we prove that
\begin{equation}\label{eq:hrnk}
	|u(t)|^2 + |u^{[1]}(t)|^2 \to 0 \qquad \mbox{ as }\ |t|\to\infty;
\end{equation}
on the other hand, \eqref{eq:hrnk} does not hold for $u_m$ by
the so-called lemma on three periods. This contradicts the closeness of $u$
and $u_m$; therefore $u$ cannot be an eigenfunction and $\la$ cannot be an
eigenvalue of $S$.

The paper is organized as follows. Some auxiliary results
(a structure theorem for the space $W^{-1}_{2,unif}(\bR)$,
relation~\eqref{eq:hrnk} for eigenfunctions, the lemma on three
periods and some estimates for solutions of equation~\eqref{eq:nprd})
are established in Section~\ref{sec:pre}. The Gronwall-type
estimates for closeness of solutions $u$ and $u_m$ of
equations~\eqref{eq:nprd} and \eqref{eq:prd} and the proof of
Theorem~\ref{thm:main} are given in Section~\ref{sec:prf}. Finally,
Theorem~\ref{thm:sgp} is proved in Section~\ref{sec:appl}.

\section{Preliminaries and some auxiliary results}\label{sec:pre}

We collect in this section some auxiliary results to be exploited
later on. In the following, $W^s_2(\bR)$, $s\in\bR$, will denote
the Sobolev spaces~\cite{Sob}, $\|\cdot\|$ without any subscript will
always stand for the $L_2(\bR)$-norm and $f^{[1]}$ for the
{\em quasi-derivative} $f' - \si f$ of a function~$f$.

\subsection{The structure of the space $W^{-1}_{2,unif}(\bR)$}\label{ssec:str}
We recall that $W^{-1}_2(\bR)$ is the dual space to the Sobolev space
$W^1_2(\bR)$, i.e., it consists of those {\em distributions}~\cite{dstr},
which define continuous functionals on $W^1_2(\bR)$.
With $\lan \cdot,\cdot\ran$ denoting the duality, we have
for $f\in W^{-1}_2(\bR)$
\[
	\|f\|_{W^{-1}_2(\bR)} := \sup_{0\ne\psi\in W^1_2(\bR)}
			\frac{|\lan \psi, f \ran|}{\|\psi\|_{W^1_2(\bR)}}.
\]
The local uniform analogue of this space is constructed as follows.
Put
\begin{equation}\label{eq:phi}
   \phi(t):=  \left\{ \begin{array}{ll}
		2(t+1)^2\quad &\mbox{if\quad $t\in [-1,-1/2)$},\\
		1-2t^2\quad   &\mbox{if\quad $t\in [-1/2,1/2)$},\\
		2(t-1)^2\quad &\mbox{if\quad $t\in [1/2,1]$},\\
	        0\quad        &\mbox{otherwise}, 	  \end{array}
              \right.
\end{equation}
and $\phi_n(t) := \phi(t-n)$ for $n\in\bZ$. We say that $f$ belongs to
$W^{-1}_{2,unif}(\bR)$ if $f\phi_n$ is in $W^{-1}_{2}(\bR)$ for all $n\in\bZ$
and
\[
	\|f\|_{W^{-1}_{2,unif}(\bR)} := \sup_{n\in\bZ}
		\|f \phi_n\|_{W^{-1}_{2}(\bR)} < \infty.
\]
The space $W^{-1}_{2,unif}(\bR)$ has the following structure (see details
in~\cite{HM}).

\begin{theorem}[{\cite{HM}}]\label{thm:str}
For any $f\in W^{-1}_{2,unif}(\bR)$ there exist functions
$\si \in L_{2,unif}(\bR)$ and $\tau\in L_{1,unif}(\bR)$ such that
$f = \si'+\tau$ and
\begin{equation}\label{eq:norms}
	C^{-1} \bigl(\|\si\|_{2,unif} + \|\tau\|_{1,unif} \bigr) \le
			\|f\|_{W^{-1}_{2,unif}(\bR)} \le
	C      \bigl(\|\si\|_{2,unif} + \|\tau\|_{1,unif} \bigr)
\end{equation}
with some constant $C$ independent of $f$.
Moreover, the function $\tau$ can be chosen uniformly bounded.
\end{theorem}

We say that a distribution $f$ is $T$-periodic if
\(
	\lan f, \psi (t) \ran = \lan f , \psi (t+T) \ran
\)
for any test function $\psi$. It follows from the proof of
Theorem~\ref{thm:str} that for a $T$-periodic
potential $f\in W^{-1}_{2,unif}(\bR)$ the functions $\si$ and $\tau$ may
be taken $T$-periodic and constant, respectively (see~\cite[Remark~2.3]{HM}).

\subsection{Inequalities, quasi-derivatives etc}
The aim of the next three lemmata is to prove that any eigenfunction vanishes
at infinity together with its quasi-derivative. For Schr\"odinger operators
with regular potentials this result is usually derived from Harnack's
inequality (see, e.g., \cite{hrnk}).

\begin{lemma}\label{lem:embed}
Suppose that $f\in W^1_2(\bR)$; then
\[
	\max_{t\in\bR} |f(t)| \le \|f\|_{W^1_2(\bR)}.
\]
\end{lemma}

\begin{proof}
Since $f'\ov f \in L_1(\bR)$, the relation
\[
	|f(t)|^2 - |f(s)|^2 = 2 \int_s^t \Re f'\ov f
\]
implies that
\(
	\lim_{s\to-\infty} |f(s)| = 0.
\)
Therefore
\[
	|f(t)|^2 = 2 \int_{-\infty}^t \Re f'\ov f \le
		\int_{-\infty}^t (|f'|^2 + |f|^2) \le \|f\|^2_{W^1_2(\bR)}
\]
and the proof is complete.
\end{proof}

\begin{lemma}\label{lem:qdr}
The quasi-derivative $f^{[1]}:= f' - \si f$ of any function $f\in\fD(S)$
belongs to $L_2(\bR)$.
\end{lemma}

\begin{proof}
It suffices to show that $\si f \in L_2(\bR)$ for any $\si\in L_{2,unif}(\bR)$
and any $f\in W^1_2(\bR)$. Let $\phi$ be defined through~\eqref{eq:phi};
then $0\le\phi\le1$, $|\phi'|\le2$, and
\(\sum_{n=-\infty}^\infty \phi(t-n) \equiv1\).
Denote $\phi_n(t) := \phi(t-n)$ and $f_n(t):= f(t)\phi_n(t)$;
then $\sum_{n\in\bZ} f_n =f$,
\[
	\|f_n\|^2_{W^1_2(\bR)}
		\le 2\|f'\phi_n\|^2 + 2\|f \phi'_n\|^2 + \|f\phi_n\|^2
		\le 9\|f\|^2_{W^1_2[n-1,n+1]},
\]
and hence
\[
	\sum_{n\in\bZ} \|f_n\|^2_{W^1_2(\bR)} \le 18 \|f\|^2_{W^1_2(\bR)}
		< \infty.
\]

Recalling now Lemma~\ref{lem:embed}, we derive the inequalities
\begin{align*}
	\|\si f\|^2 &
		\le 2\Bigl(\sum_{n\in\bZ} \|\si f_{2n}\|\Bigr)^2
		  + 2\Bigl(\sum_{n\in\bZ} \|\si f_{2n-1}\|\Bigr)^2 \\
	& \le 2 \sum_{n\in\bZ} \|\si f_{n}\|^2
	\le 2\sum_{n\in\bZ}\|\si\|^2_{L_2(n-1,n+1)}\max|f_n(t)|^2\\
	& \le 2\sum_{n\in\bZ}\|\si\|^2_{L_2(n-1,n+1)}\|f_n\|^2_{W^1_2(\bR)}
	\le 4\|\si\|^2_{2,unif} \sum_{n\in\bZ}\|f_n\|^2_{W^1_2(\bR)} <\infty,
\end{align*}
whence $\si f \in L_2(\bR)$. The lemma is proved.
\end{proof}

\begin{lemma}\label{lem:ef}
Suppose that $u$ is an eigenfunction of the \op~$S$ corresponding to
an \ev~$\la$. Then
\[
	|u(t)|^2 + |u^{[1]}(t)|^2 \to 0 \quad \mbox{as}\quad |t|\to\infty.
\]
\end{lemma}

\begin{proof}
Since $u \in W^1_2(\bR)$, we have $\lim_{|t|\to\infty}u(t)=0$ and
\[
	\int_t^{t+1}|\si u'| \le
		\Bigl(\int_t^{t+1}|\si|^2 \Bigr)^{1/2}
		\Bigl(\int_t^{t+1}|u'|^2 \Bigr)^{1/2}
	    \le \|\si\|_{2,unif} \Bigl(\int_t^{t+1}|u'|^2 \Bigr)^{1/2} \to 0
\]
as $|t|\to\infty$. Hence
\begin{align*}
	\sup_{\al,\be\in[t,t+1)} |u^{[1]}(\be) - u^{[1]}(\al)|
		&\le \int_t^{t+1} |(u^{[1]})'| \\
		&\le \int_t^{t+1} |\si u'| + \int_t^{t+1} |\tau u|
			+ \int_t^{t+1} |\la u| \to 0
\end{align*}
as $|t|\to\infty$, which implies existence of the limits
\(
	u_\pm := \lim_{t\to\pm\infty} u^{[1]}(t).
\)
By Lemma~\ref{lem:qdr} $u^{[1]}\in L_2(\bR)$, whence $u_\pm=0$ and the proof
is complete.
\end{proof}

\subsection{Some results on first order differential systems}
It is convenient to rewrite the eigenvalue equation~\eqref{eq:nprd}
as a first order system
\begin{equation}\label{eq:syst}
	\frac{d}{dt}\bn{u_1}{u_2} =
		\begin{pmatrix} \si & 1 \\ -\si^2 + \tau -\la & -\si
		\end{pmatrix} \bn{u_1}{u_2}.
\end{equation}
A function $u$ is easily seen to be a solution
of~\eqref{eq:nprd} if and only if the vector-function $U:= (u,u^{[1]})^T$
solves~\eqref{eq:syst}. Lemma~\ref{lem:ef} states that if $u$ is in addition
an eigenfunction of~$S$, then $|U(t)| \to 0$ as $|t|\to\infty$,
where for a vector $X\in\bC^2$ we denote by $|X|$ its Euclidean norm.
The following statement (known as the ``lemma on three periods", cf.~\cite{G,DS})
shows that for periodic potentials $q =\si'+\tau\in W^{-1}_{2,unif}(\bR)$ no nontrivial solution
of~\eqref{eq:syst} vanishes at infinity (and hence by Lemma~\ref{lem:ef} the
corresponding Schr\"odinger operators have empty point spectrum; in fact, the
spectrum is then absolutely continuous~\cite{HM}).

\begin{lemma}\label{lem:3per}
Suppose that $A$ is a $T$-periodic $2\times2$ matrix with locally integrable
entries and zero trace and let $U$ be any solution of the equation
\(
	\dot{U} = A U.
\)
Then
\[
	\max\{|U(-T)|, |U(T)|, |U(2T)|\} \ge \tfrac12 |U(0)|.
\]
\end{lemma}

\begin{proof}
Denote by $M(t)$ the monodromy matrix of the equation $\dot{U} = AU$;
then we conclude by the Liouville and Cayley-Hamilton theorems that
$\det M\equiv1$ and
\[
	M^2 + (\tr M) M + I = 0.
\]
Since $(M U)(t) = U(t+T)$, this yields
\[
	U(t+2T) + (\tr M) U(t+T) + U(t) \equiv 0.
\]
If $|\tr M| \le 1$, we take $t=0$ in the above equality to derive that
\[
	\max\{|U(2T)|, |U(T)|\} \ge \tfrac12|U(0)|,
\]
while otherwise $t = -T$ produces the estimate
\[
	\max\{|U(T)|, |U(-T)| \} \ge \tfrac12|U(0)|.
\]
The lemma is proved.
\end{proof}

Next we adopt the arguments from~\cite{YS} to system~\eqref{eq:syst}
to derive some {\em a priori} estimates for its solutions.

\begin{lemma}\label{lem:exp}
Let $a$ and $b$ be some locally summable functions and $X(t)$ be an
arbitrary nonzero solution of the equation
\[
	\frac{dX}{dt} = \begin{pmatrix} a & 1 \\ b & -a
			\end{pmatrix} X.
\]
Then for any constant $c\ge1$ and any $t\in\bR$ the following
inequality holds:
\begin{equation}\label{eq:exp}
	|X(t)| \le c \exp \Bigl\{\frac12\int_{\min\{0,t\}}^{\max\{0,t\}}
		\sqrt{ 4a(s)^2 + (c + b(s)/c)^2}\,ds \Bigr \}|X(0)|.
\end{equation}
\end{lemma}

\begin{proof}
It suffices to consider the case $t\ge 0$. Denote
\[
	A := \begin{pmatrix} a & 1 \\ b & -a
			\end{pmatrix}, \quad
	G := \begin{pmatrix} c & 0 \\ 0 & 1/c
			\end{pmatrix},\quad
	F := GA + A^* G =
		\begin{pmatrix} 2ca & c + b/c \\ c+ b/c & -2a/c
			\end{pmatrix}.
\]
Then $\xi(t):=\bigl(GX(t),X(t)\bigr)$ satisfies the relation
\[
    \dot{\xi}(t) = \bigl( G\dot{X}(t),X(t) \bigr) +
		\bigl( GX(t),\dot{X}(t) \bigr) = \bigl( FX(t),X(t) \bigr).
\]
If $\la^+$ is larger of two eigenvalues of the pencil $F-\la G$, then
$F \le \la^+ G$, whence $\dot{\xi}(t) \le \la^+(t) \xi(t)$ and
$\log \xi(t) - \log\xi(0) \le \int_0^t \la^+(s)\,ds$.
Taking the exponents of both sides, we arrive at the inequality
\[
	|X(t)|^2 \le c \xi (t)
		\le c\exp \Bigl\{ \int_0^t \la^+(s)\,ds\Bigr\} \xi(0)
		\le c^2\exp \Bigl\{ \int_0^t \la^+(s)\,ds\Bigr\} |X(0)|.
\]
To calculate $\la^+$, we observe that it is larger of two eigenvalues of
the matrix
\[
	F_1 := G^{-1/2}FG^{-1/2} =
		\begin{pmatrix} 2a & c + b/c \\ c+ b/c & -2a
			\end{pmatrix}.
\]
Since $F_1$ has a zero trace we have
\[
	\la^+ = \sqrt{-\det F_1} = \sqrt{4a^2 + (c + b/c)^2},
\]
and the proof is complete.
\end{proof}

For equation~\eqref{eq:syst}, Lemma~\ref{lem:exp} yields the following
result.

\begin{lemma}\label{lem:sol}
Suppose that the functions $\si$ and $\tau$ belong to $L_{2,loc}(\bR)$ and
$L_{1,loc}(\bR)$ respectively. Then any solution $X(t)$ of
equation~\eqref{eq:syst} satisfies the inequality
\[
   |X(t)| \le C_1  \exp \Bigl\{ \frac12\int_{\min\{0,t\}}^{\max\{0,t\}}
		(2-\eps+2\sqrt{(-\la)_+} + \si^2 + |\tau|)\,ds\Bigr\}|X(0)|,
\]
where $(-\la)_+ = \max\{-\la,0\}$ and the constants $C_1 >0$ and $\eps >0$
only depend on~$\la$.
\end{lemma}

\begin{proof}
We shall show that for any $\la\in\bR$ there exist constants
$C_1 \ge1$ and $\eps >0$ such that
\[
	4\si^2 + (C_1^2 - \si^2 + \tau - \la)^2/C_1^2
		< (2 -\eps + 2\sqrt{(-\la)_+} + \si^2 + |\tau|)^2;
\]
then the statement will follow from Lemma~\ref{lem:exp}.

For $\la \ge 0$ we put $C_1 = \sqrt{\la + 1/4} + 1/2$; then
$C_1^2 - \la = C_1$, whence
\begin{align*}
	4\si^2 + (C_1^2 - \si^2 + \tau - \la)^2/C_1^2
		&= 4\si^2 + (C_1 - \si^2 + \tau)^2/C_1^2 \\
		&\le (\si^2 + |\tau|)^2 + 2 (\si^2 + |\tau|)(2-1/{C_1})
				+ (2-1/{C_1})^2\\
		&= (\si^2 + |\tau| + 2 - 1/{C_1})^2.
\end{align*}
For $\la \in(-1,0)$ we take $C_1 = 1$; then
\begin{align*}
	4\si^2 + (1 - \si^2 + \tau - \la)^2
		&= (\si^2 - \tau)^2 -2 (\si^2 - \tau)(1-\la) +
			4\si^2 + (1-\la)^2\\
		&\le (\si^2 + |\tau| + 1- \la)^2.
\end{align*}
Finally, for $\la = -\nu^2 \le -1$ we put $C_1 = \nu = \sqrt{(-\la)_+}$; then
\begin{align*}
	4\si^2 + (C_1^2 - \si^2 + \tau - \la)^2/C_1^2
		&= 4\si^2 + (2\nu^2 - \si^2 + \tau)^2/\nu^2\\
		&= (\si^2 -\tau)^2/\nu^2 - 4(\si^2 - \tau) + 4\si^2 + 4\nu^2\\
		& \le (\si^2 +|\tau| + 2\nu)^2.
\end{align*}
The proof is complete.
\end{proof}

\begin{corollary}\label{cor:sol}
If the functions $\si$ and $\tau$ belong to $L_{2,unif}(\bR)$ and
$L_{1,unif}(\bR)$ respectively, then any solution $X(t)$ of
equation~\eqref{eq:syst} satisfies the inequality
\[
   |X(t)| \le C_1  \exp \Bigl\{(|t|+1)
	\bigl( 1 - \eps + \sqrt{(-\la)_+} + \tfrac12\|\si\|_{2,unif}^2 +
		\tfrac12\|\tau\|_{1,unif}\bigr)
			\Bigr\}|X(0)|
\]
with suitable $C_1>0$ and $\eps>0$ independent of $\si$, $\tau$, and $t$.
\end{corollary}

\section{Proof of Theorem~\ref{thm:main}}%
\label{sec:prf}

Denote by $\ga_0$ the lower bound of the \op~$S$ and by $\ga$ any number
less than~$\ga_0$. We fix now any $\la > \ga$ and shall prove
that $\la$ is not an eigenvalue of~$S$. The eigenvalue equation
$ l(u) = \la u $ can be recast as a first order system
\[
	\frac{d}{dt}\bn{u_1}{u_2} =
		\begin{pmatrix} \si & 1 \\ -\si^2 + \tau - \la & -\si
		\end{pmatrix} \bn{u_1}{u_2}
\]
with $u_1 = u$ and $u_2 = u^{[1]}$. Set
\[
	U(t) = \bn{u_1}{u_2}, \qquad
	A = \begin{pmatrix} 0 & 1 \\ \ga-\la & 0 \end{pmatrix},\qquad
	B(t) = \begin{pmatrix} \si & 0 \\ -\si^2+\tau-\ga & -\si\end{pmatrix};
\]
then
\[
	\frac{dU}{dt} = AU + BU.
\]
Let $\ti U$ be a solution of the equation
\[
	\frac{d\ti U}{dt} = A{\ti U} + \ti B \ti U
\]
with the initial condition $\ti U(0) = U(0)$; here
\[
	\ti U(t) = \bn{\ti u_1}{\ti u_2}, \qquad
	\ti B(t) = \begin{pmatrix}
		\ti \si & 0 \\ -\ti\si^2+\ti\tau-\ga & -\ti\si
		\end{pmatrix}.
\]
In the following $|U(t)|$, $|A|$, $|B(t)|$ etc. will denote the Euclidean norms
of the vector $U(t)$ and matrices $A$ and $B(t)$ respectively.

\begin{lemma}\label{lem:Gronw}
There exists a constant $C_2=C_2(\la)>0$ such that for all $t\in\bR$ the
following inequality holds:
\begin{equation}\label{eq:Gronw}
	|U(t) - \ti U(t)| \le C_2
	\int_{\min\{0,t\}}^{\max\{0,t\}} |B(s) - \ti B(s)| |\ti U(s)|\,ds\
	    \exp\Bigl\{\int_{\min\{0,t\}}^{\max\{0,t\}} |B(s)|\,ds \Bigr\}.
\end{equation}
\end{lemma}

\begin{proof} To be definite, we consider the case $t>0$.
The function $V:= U -\ti U$ solves the equation
\[
	\frac{dV}{dt} = AV + BV + (B- \ti B)\ti U
\]
and satisfies the initial condition $V(0)=0$, whence
\[
	V(t) = \int_0^t e^{(t-s)A}B(s)V(s)\,ds +
		\int_0^t e^{(t-s)A}\bigl(B(s) - \ti B(s)\bigr)\ti U(s)\,ds.
\]
Since $\la - \ga > 0$, the group $e^{tA}$ is uniformly bounded by some
constant $C_3$, whence
\[
	|V(t)| \le C_3 \int_0^t |B(s)||V(s)|\,ds +
		C_3\int_0^t |B(s)-\ti B(s)||\ti U(s)|\,ds.
\]
The Gronwall inequality now yields
\[
	|V(t)| \le C_2 \int_0^t |B(s) - \ti B(s)| |\ti U(s)|\,ds\
			\exp\Bigl\{\int_0^t |B(s)|\,ds \Bigr\}
\]
with $C_2 = C_3 e^{C_3}$, and the lemma is proved.
\end{proof}

\begin{lemma}\label{lem:norm}
For any $a,b\in\bC$,
\[
	\left|\begin{pmatrix} a & 0 \\ b & -a \end{pmatrix} \right|
	\le |b|/2 + \bigl(|b|^2/4 + |a|^2 \bigr)^{1/2}.
\]
\end{lemma}

\begin{proof}
We have to maximize the expression
\[
	\left|\begin{pmatrix} a & 0 \\ b & -a \end{pmatrix}
					\bn{z_1}{z_2}\right|^2 =
		|az_1|^2 + |bz_1|^2 - 2 \Re a\ov b z_1 \ov z_2 + |az_2|^2
\]
over $(z_1,z_2)^T \in \bC^2$ of unit length. Introducing $\th \in [0,\pi/2]$
via $|z_1| = \cos \th$, $|z_2| = \sin \th$ and using the relations
$2|z_1z_2| = \sin 2\th$, $|z_1|^2 = (1 + \cos2\th)/2$,
we bound the above expression by
\begin{align*}
	|a|^2 + |b|^2(1 + \cos 2\th)/2 + |a||b|\sin2\th
	&\le |a|^2 +  |b|^2/{2} +
		|b| \bigl(|b|^2/4 + |a|^2\bigr)^{1/2} \\
	&= \left(|b|/2 + \bigl(|b|^2/4 + |a|^2 \bigr)^{1/2} \right)^2,
\end{align*}
and the proof is complete.
\end{proof}

\begin{proofof}{Proof of Theorem~\ref{thm:main}.}
Assume that the assumptions of the theorem hold and $\la>\ga$ is an
eigenvalue of the \op~$S$ with the corresponding eigenfunction $u$
normalized by the condition $|u(0)|^2 + |u^{[1]}(0)|^2=1$.
Fix a sequence $q_m = \si'_m + \tau_m$ of $T_m$-periodic potentials
from the class $W^{-1}_{2,unif}(\bR)$ such that $T_m\to\infty$ and
\begin{equation}\label{eq:prox}
   \lim_{m\to\infty} \exp (C T_m)
 	\left\{ \Bigl(\int_{-T_m}^{2T_m} |\si - \si_m|^2\Bigr)^{1/2} +
		\int_{-T_m}^{2T_m}|\tau - \tau_m|
	\right\} = 0
\end{equation}
for any $C < \infty$.

Put $U = (u,u^{[1]})^T$ and denote by $U_m$ a solution of equation
\[
	\frac{dU_m}{dt} = (A + B_m) U_m, \qquad
		B_m(t) := \begin{pmatrix}
			\si_m & 0 \\ -\si_m^2 + \tau_m -\ga & -\si_m
			\end{pmatrix},
\]
with the initial condition $U_m(0) = U(0)$; then by Lemma~\ref{lem:3per}
\begin{equation}\label{eq:below}
	\max\{|U_m(-T_m)|,|U_m(T_m)|,|U_m(2T_m)|\}\ge 1/2.
\end{equation}
According to Corollary~\ref{cor:sol}
\[
   |U_m(t)| \le C_1 \exp\bigl\{(|t|+1)(1+|\ga|^{1/2} +
	\tfrac12\|\si_m\|^2_{2,unif} + \tfrac12\|\tau_m\|_{1,unif}) \bigr\}
\]
and therefore~\eqref{eq:Gronw} gives
\begin{multline}
    |U(t)-U_m(t)| \le C_1 C_2
	\int_{\min\{0,t\}}^{\max\{0,t\}}|B(s) - B_m(s)|\,ds\
	\exp\Bigl\{ \int_{\min\{0,t\}}^{\max\{0,t\}} |B(s)|\,ds \Bigr\}\times\\
	\times\exp\Bigl\{(|t|+1)(1+|\ga|^{1/2} +
	\tfrac12\|\si_m\|^2_{2,unif} + \tfrac12\|\tau_m\|_{1,unif})\Bigr\}.
    \label{eq:Udif}
\end{multline}

The integrals in~\eqref{eq:Udif} can be estimated by means of
Lemma~\ref{lem:norm} as follows. First, we see that
\[
	|B-B_m| \le |\tau-\tau_m| + |\si - \si_m|(|\si| + |\si_m|+1),
\]
whence
\begin{multline}\label{eq:Bdif}
	\int_{\min\{0,t\}}^{\max\{0,t\}} |B(s)-B_m(s)|\,ds \le
	\int_{\min\{0,t\}}^{\max\{0,t\}}|\tau-\tau_m| \\
	+ \Bigl(\int_{\min\{0,t\}}^{\max\{0,t\}}
		3(|\si|^2 + |\si_m|^2 + 1)\Bigr)^{1/2}
	  \Bigl(\int_{\min\{0,t\}}^{\max\{0,t\}} |\si - \si_m|^2\Bigr)^{1/2}.
\end{multline}
In the same manner we conclude that
\[
	|B| \le 1 + |\ga| + |\si|^2 + |\tau|,
\]
which implies
\begin{equation}\label{eq:B}
	\int_{\min\{0,t\}}^{\max\{0,t\}} |B(s)|\,ds \le (|t|+1)
		(1 + |\ga| + \|\si\|^2_{2,unif} + \|\tau\|_{1,unif}).
\end{equation}

Since $\|\si_m\|^2_{2,unif} \le 2\|\si\|^2_{2,unif}$ and
$\|\tau_m\|_{1,unif} \le 2\|\tau\|_{1,unif}$ for all $m$ large enough,
combination of relations \eqref{eq:Udif}--\eqref{eq:B}
yields the inequality
\begin{multline*}
	|U(t) - U_m(t)| \le
		C_3\exp \bigl\{|t| (2 + |\ga|^{1/2} + |\ga|
			+ 2\|\si\|^2_{2,unif} + 2\|\tau\|_{1,unif}) \bigr\}\times\\
	\times\Bigl\{ \Bigl(\int_{\min\{0,t\}}^{\max\{0,t\}}
				 |\si - \si_m|^2\Bigr)^{1/2} +
		\int_{\min\{0,t\}}^{\max\{0,t\}} |\tau - \tau_m|
	\Bigr\}
\end{multline*}
with some $C_3$ independent of $m$. Substituting now
$t = \pm T_m,\, 2T_m$ and recalling~\eqref{eq:prox}, we find that
\[
	\lim_{m\to\infty} \max_{t = \pm T_m, 2T_m} |U(t) - U_m(t)| =0.
\]
By virtue of~\eqref{eq:below} this implies that
$U(t)$ does not tend to $0$ as $|t|\to\infty$.
Therefore $u$ is not an eigenfunction of~$S$ and hence
the point spectrum of $S$ is empty. The theorem is proved.
\end{proofof}

\begin{remark}\label{rem:q} \rm
Observe that we only need \eqref{eq:prox} to hold for
\[
   C < C_q :=
       4 + 2|\ga|^{1/2} + 2|\ga| + 4\|\si\|^2_{2,unif} + 4\|\tau\|_{1,unif}.
\]
Recall~\cite{HM} that $\ga_0 \ge - (a\|q\|_{W^{-1}_{2,unif}(\bR)} + b)^4$
with certain $a,b>0$ independent of $q$ and that
$\|\si\|^2_{2,unif} \le 64\|q\|_{W^{-1}_{2,unif}(\bR)}$
and $\|\tau\|_{1,unif} \le 3\|q\|_{W^{-1}_{2,unif}(\bR)}$;
therefore for some $a',b'>0$ independent of $q$ we have
\[
	C_q \le (a'\|q\|_{W^{-1}_{2,unif}(\bR)} + b')^4.
\]
\end{remark}

\section{Application to quasiperiodic potentials}\label{sec:appl}

In this section, we shall establish Theorem~\ref{thm:sgp} showing which
of the functions~\eqref{eq:qper} are singular Gordon potentials.
We shall prove first some auxiliary results.

\begin{lemma}\label{lem:add}
Suppose that $f\in W^1_2(a,b)$ with $b-a\ge 1$; then
for any $c\in (a,b)$ and any $\eps\in(0,b-c)$ the following inequality holds:
\[
	\int_a^c |f(t+\eps) - f(t)|^2\,dt \le 7 \eps^2 \|f\|^2_{W^1_2(a,b)}.
\]
\end{lemma}

\begin{proof}
We put
\[
   \psi(t)=    \left\{ \begin{array}{ll}
	f(t)\quad &\mbox{if\quad $t\in [a,b]$},\\
	f(2a-t)\,\dfrac{t+b-2a}{b-a} \quad &\mbox{if\quad $t\in (2a-b,a)$},\\
	f(2b-t)\,\dfrac{2b-a-t}{b-a} \quad &\mbox{if\quad $t\in (b,2b-a)$},\\
	0 \qquad   &\mbox{otherwise.}
               		\end{array} \right.
\]
Then $\psi \in W^1_2(\bR)$ and
\[
	\|\psi\|^2_{W^1_2(\bR)} =
		\|\psi\|^2_{L_2(\bR)} + \|\psi'\|^2_{L_2(\bR)} \le
		7 \left(\|f\|^2_{L_2(a,b)}+\|f'\|^2_{L_2(a,b)}\right)
		= 7 \|f\|^2_{W^1_2(a,b)}.
\]

Denote by $\hat \psi$ the Fourier transform of $\psi$. Then using Plancherel's
theorem and recalling the equivalent definition of the norm in $W^1_2(\bR)$,
we get
\begin{align*}
 	\int_a^c |f(t+\eps) - f(t)|^2\,dt &\le
		\int_\bR |\psi(t+\eps) - \psi(t)|^2\,dt =
		\int_\bR |\hat \psi(u)|^2|e^{i\eps u}-1|^2\,du\\
	& \le \eps^{2} \max_{u \in \bR}\frac{|e^{i\eps u}-1|^2}{|\eps u|^{2}}
	    \int_\bR (1+u^2) |\hat \psi(u)|^2\,du
			\le  \eps^{2} 	\|\psi\|^2_{W^1_2(\bR)},
\end{align*}
and the lemma follows.
\end{proof}

\begin{lemma}\label{lem:mult}
Suppose that $f \in W^1_2(\bR)$ and $a>1$; then for any $b\ge e$
\[
	\int_1^b|f(t) - f(at)|^2\,dt
		\le 7 a b^2(a-1)^{2}\|f\|^2_{W^1_2(\bR)}.
\]
\end{lemma}

\begin{proof}
Put $c= \log ab$ and $g(t):=f(e^t)$; then $g\in W^1_2(0,c)$ and
by Lemma~\ref{lem:add} we have
\begin{align*}
	\int_1^b |f(t) - f(at)|^2\,dt
	&= \int_0^{\log b} |g(u) - g(u+\log a)|^2e^u\,du \\
	& \le b \int_0^{\log b} |g(u) - g(u+\log a)|^2 \,du \le
	7 b (\log a)^{2} \|g\|^2_{W^1_2(0,c)}.
\end{align*}
Observe that
\[
	\int_0^c|g(u)|^2\,du = \int_0^{\log ab}|f(e^u)|^2\,du
	= \int_1^{ab}|f(t)|^2\frac{dt}{t} \le \|f\|^2_{L_2(\bR)}
\]
and
\[
	\int_0^c|g'(u)|^2\,du = \int_0^{\log ab}|f'(e^u)|^2e^{2u}\,du
	= \int_1^{ab}|f'(t)|^2 t\,dt \le ab \, \|f'\|^2_{L_2(\bR)}.
\]
Therefore
\[
	\|g\|^2_{W^1_2(0,c)} \le ab\, \|f\|^2_{W^1_2(\bR)}
\]
and the result follows.
\end{proof}

\begin{lemma}\label{lem:albe}
For any $s \in [0,1]$ there exists a constant $C_s>0$ such that the inequality
\[
	\int_{-T}^{2T} |f(\al t + \th) - f(\be t + \th)|^2\,dt
		\le C_s T^{2s}\al^{-1}(\be-\al)^{2s} \|f\|^2_{W^s_2(\bR)}.
\]
holds for all $f \in W^s_2(\bR)$, all
$\al,\be >0$, $\al \le \be \le 2\al$, all $\th\in\bR$, and all $T\ge 1/\al$.
\end{lemma}

\begin{proof}
Denote by $A_s$ an operator from $W^s_2(\bR)$ into $L_2(-T,2T)$ defined by
\[
	(A_sf) (t) := f(\al t +\th) - f(\be t + \th);
\]
then the statement of the lemma asserts that
\begin{equation}\label{eq:interp}
	\|A_s\|^2 \le C_s T^{2s}\al^{-1}(\be-\al)^{2s}
\end{equation}
with some constant $C_s>0$ independent of $\al,\be$,$\th$, and $T$.
We shall prove inequality~\eqref{eq:interp} for $s=0$ and $s=1$ and then
interpolate (see~\cite[Ch.1]{LM} for details) to cover all $s$ in between.

It is easily seen that
\[
	\|A_0f\|^2_{L_2(-T,2T)} \le
		2\frac{\al+\be}{\al\be}\|f(t)\|^2_{L_2(\bR)}
		\le 4\al^{-1} \|f(t)\|^2_{L_2(\bR)},
\]
whence~\eqref{eq:interp} holds for $s=0$ with $C_0 = 4$.

Suppose now that $f \in W^1_2(\bR)$ and put
$g(t) = f(\al T t - 2\be T + \th)$; then $g \in W^1_2(\bR)$ and
\[
	\|g\|^2_{W^1_2(\bR)} = \frac1{\al T}\|f\|^2_{L_2(\bR)} +
		\al T \|f'\|^2_{L_2(\bR)} \le \al T \|f\|^2_{W^1_2(\bR)}.
\]
Moreover,
\begin{align*}
	\|A_1f\|^2_{L_2(-T,2T)} &=
	\int_{-T}^{2T} |f(\al t + \th) - f(\be t + \th)|^2\,dt \\
	&= T\int_1^4 \bigl|g\bigl(u + 2\tfrac{\be - \al}{\al}\bigr) -
		g \bigl(\tfrac{\be}{\al}u\bigr)\bigr|^2\,du \\
	& \le 2T \int_1^4 \bigl|g\bigl(u + 2{\tfrac{\be - \al}{\al}}\bigr)
		-g(u)\bigr|^2\,du +
	2T \int_1^4 \bigl|g(u)- g\bigl(\tfrac{\be}{\al}u\bigr) \bigr|^2\,du.
\end{align*}
These two summands are estimated by Lemmata~\ref{lem:add} and \ref{lem:mult}
as
\begin{align*}
    2T \int_1^4 \bigl|g\bigl(u + 2\tfrac{\be - \al}{\al}\bigr)
		-g(u)\bigr|^2\,du
      &\le 14T \bigl(2\tfrac{\be - \al}{\al} \bigr)^{2}\|g\|^2_{W^1_2(\bR)},\\
    2T \int_1^4 \bigl|g(u)- g\bigl(\tfrac{\be}{\al}u\bigr) \bigr|^2\,du
	&\le 14T\,4^2\tfrac{\be}{\al}
	\bigl(\tfrac{\be - \al}{\al} \bigr)^{2}\|g\|^2_{W^1_2(\bR)}
\end{align*}
respectively. Combining the above inequalities we finally get
\[
	\|A_1f\|^2_{L_2(-T,2T)} \le 14\cdot36 \,T^2 \al^{-1} (\be - \al)^2
		\|f\|^2_{W^1_2(\bR)},
\]
and~\eqref{eq:interp} holds for $s=1$ with $C_1 = 14 \cdot 36$.
It suffices now to put $C_s = C_0^{1-s}C_1^s$, and the proof is complete.
\end{proof}

We denote by $W^s_{2,unif}(\bR)$ the set of all functions
$f\in L_{2,loc}(\bR)$ such that $f \phi_n \in W^s_2(\bR)$ for any $n\in \bZ$
and
\[
	\|f\|_{W^s_{2,unif}(\bR)} :=
		\sup_{n\in\bZ}\|f \phi_n\|_{W^s_2(\bR)} < \infty;
\]
here $\phi_n(t) = \phi (t-n)$ and $\phi$ is the function of~\eqref{eq:phi}.

\begin{theorem}\label{thm:hldr}
Suppose that $f\in W^s_{2,unif}(\bR)$, $s\in[0,1]$, $\th \in \bR$,
and $\al,\be>0$, $\al \le \be \le 2\al$. Then for any
$T\ge \max\{ 1/\al, |\th|\}$ we have
\[
	\int_{-T}^{2T} |f(\al t + \th) - f(\be t + \th)|^2\,dt
		\le C T^{2s+1} (\be-\al)^{2s} \|f\|^2_{W^s_{2,unif}(\bR)},
\]
where $C>0$ is a constant independent of $f$, $\al$, $\be$, and $T$.
\end{theorem}

\begin{proof}
Put
\(
	g = f \sum_{k=-n}^{2n} \phi_n
\)
with $n = [\be T]+2$, where $[a]$ is the integral part of a number $a$.
Then $g\in W^s_2(\bR)$ and
\[
	\|g\|_{W^s_2(\bR)} \le (3n+1)\|f\|_{W^s_{2,unif}(\bR)}
		\le 10 \be T \|f\|_{W^s_{2,unif}(\bR)}.
\]
Moreover, $f(t) = g(t)$ for $t\in [-\be T + \th,2\be T + \th]$ provided
$|\th| \le T$, so we then have
\[
	\int_{-T}^{2T} |f(\al t + \th) - f(\be t + \th)|^2\,dt =
	\int_{-T}^{2T} |g(\al t + \th) - g(\be t + \th)|^2\,dt,
\]
and the claim follows from Lemma~\ref{lem:albe}.
\end{proof}

We recall that an irrational number $\al$ is a Liouville number~\cite{Lio} if
there exists a sequence of rationals $\al_m = R_m/T_m$ such that
\begin{equation}\label{eq:lio}
	|\al - \al_m| \le C m^{-T_m}
\end{equation}
with a suitable constant $C$.
Combining the above results, we can now prove Theorem~\ref{thm:sgp}.

\begin{proofof}{Proof of Theorem~\ref{thm:sgp}. }
Suppose that
\[ 	
    q(t) = \si'_1(t) + \si'_2(\al t + \th) + \tau_1(t) + \tau_2(\al t + \th),
\]	
where $\si_1,\si_2 \in L_{2,loc}(\bR)$ and $\tau_1, \tau_2\in L_{1,loc}(\bR)$
are $1$-periodic,  $\al,\th \in [0,1)$ and $\al$ is a Liouville number.

First of all we represent the function $\tau_2$ as $\si'_3 + c$
where $c \equiv \int_0^1 \tau_2$ and $\si_3$ is the $1$-periodic primitive
of $\tau_2 - c$. Then $\si_3 \in W^1_{1,loc}(\bR)$ and hence
$\si_3 \in W^s_{2,loc}(\bR)$ with any $s<1/2$ by the Sobolev embedding
theorem~\cite{Sob}. Therefore
\[
	q(t) = \si'_1(t) + \ti\si'_2(\al t + \th) + \ti\tau_1 (t),
\]
where $\ti\si_2 = \si_2 + \si_3$ and $\ti\tau_1 = \tau_1 + c$.
We approximate $q$ by a sequence of $T_m$-periodic potentials
\[
	q_m(t) := \si'_1(t) + \ti\si'_2(\al_m t + \th) + \ti\tau_1 (t)
\]
where $\al_m = R_m/T_m$ is the $m$-th approximate of $\al$
satisfying~\eqref{eq:lio}. Combination of
Theorem~\ref{thm:hldr} and relation~\eqref{eq:lio} shows that $q$ is a
singular Gordon potential, and the proof is complete.
\end{proofof}


\end{document}